\documentclass[11pt,reqno]{amsart}
\usepackage[margin=2.5cm]{geometry}

\usepackage{graphicx} 
\usepackage{amsmath}
\usepackage{amssymb}
\usepackage{amsthm}
\usepackage[foot]{amsaddr}
\usepackage{mathtools}
\usepackage{float}
\usepackage{tikz}
\usetikzlibrary{shapes.geometric, arrows, positioning}
\usepackage{hyperref, cleveref}

\newtheorem{remark}{Remark}

\usepackage{bbm}
\usepackage{subcaption}

\newcommand\numberthis{\addtocounter{equation}{1}\tag{\theequation}}
\usepackage{algorithm}

\usepackage{xcolor}

\usepackage[backend=bibtex, bibstyle=ieee]{biblatex}
\usepackage{mathtools, stmaryrd}
\usepackage{xparse} \DeclarePairedDelimiterX{\Iintv}[1]{\llbracket}{\rrbracket}{\iintvargs{#1}}
\NewDocumentCommand{\iintvargs}{>{\SplitArgument{1}{,}}m}
{\iintvargsaux#1} %
\NewDocumentCommand{\iintvargsaux}{mm} {#1\mkern1.5mu,\mkern1.5mu#2}

\def\bpi{\boldsymbol{\pi}}

\title[Battery Management Optimization]{Maximizing On-Bill Savings through Battery Management Optimization}
\author{Rene Carmona$^1$, Xinshuo Yang$^1$, Siddharth Bhela$^2$}
\address{$^1$ Princeton University, Princeton, USA}
\author{Claire Zeng$^1$}
\address{$^2$ Siemens Technology, Princeton, NJ, USA}
\email{\{rcarmona,xy3134\}@princeton.edu, siddharth.bhela@siemens.com, cszeng@princeton.edu}
\date{September 2024}

\begin{document}

\maketitle
\def\thefootnote{$^1$}\footnotetext{These three authors contributed equally to this work.}\def\thefootnote{\arabic{footnote}}

\begin{abstract}
    In many power grids, a large portion of the energy costs for commercial and industrial consumers are set with reference to the coincident peak load, the demand during the maximum system-wide peak, and their own maximum peak load, the non-coincident peak load. Coincident-peak based charges reflect the allocation of infrastructure updates to end-users for increased capacity, the amount the grid can handle, and for improvement of the transmission, the ability to transport energy across the network. Demand charges penalize the stress on the grid caused by each consumer's peak demand. Microgrids with a local generator, controllable loads, and/or a battery technology have the flexibility to cut their peak load contributions and thereby significantly reduce these charges. This paper investigates the optimal planning of microgrid technology for electricity bill reduction. The specificity of our approach is the leveraging of a scenario generator engine to incorporate probability estimates of coincident peaks and non-coincident peaks into the optimization problem. 
\end{abstract}

\section{Introduction} 

Electric utilities provide various services related to the generation, transmission, and distribution of electricity. Utilities are also responsible for accurate metering and billing of customers based on energy consumption. Bills for commercial and industrial customers are typically generated monthly and consist of various charges \cite{PSEG}: i) an energy charge for the total electricity consumed based on a fixed cents/kWh rate; ii) a non-coincident peak (NCP) demand charge based on the maximum power consumed in a given $15$ or $30$-minute interval during the billing cycle (depending on the utility, this demand can be the highest across multiple months or billing cycles); iii) a delivery charge that covers the cost of delivering electricity and costs for maintenance of infrastructure like power lines, substations, and transformers; and iv) a capacity or coincident peak (CP) charge or discount that is based on the usage during the regional transmission operator (RTO), independent system operators (ISO), and/or local utility’s highest net demand periods - these are typically the top one or more peak hours of the year or season. For instance, Public Service Enterprise Group (PSEG) charges are based on the single highest net demand period referred to as the 1CP \cite{PSEG}, The Pennsylvania-New Jersey-Maryland Interconnection (PJM) charges are based on the $5$ highest hours (5CP) \cite{PJM5CP}, and Electric Reliability Council of Texas (ERCOT) is based on the $4$ highest hours (4CP) \cite{ERCOT}. Note that both the NCP and CP charges are determined ex-post, i.e., the periods during which the highest local or system-time demand occurs are unknown beforehand. However, they can be forecasted \cite{carmona2024coincidentpeakpredictioncapacity}, \cite{DOW18}. It is important to note that the NCP and CP periods may not coincide. In fact, even the CP period for the local utility may not coincide with the RTO/ISO CP \cite{carmona2024coincidentpeakpredictioncapacity}. Utility bills can also have other fixed costs such as system benefits charges, societal benefit charges, and regulatory charges. 

Of all the billing charges and discounts, the CP charges can represent a significant portion of the electricity bill. The exact percentage varies depending on several factors, including the customer's usage patterns, the specific rates of their utility or transmission company, and their consumption during the identified peak periods. Generally, CP charges can constitute anywhere from 20\% to 70\% of a large customer's total annual electricity costs \cite{LIU13}. With the increasing penetration of intermittent renewable energy sources and the demand growth from data centers, these costs are expected to increase even more. Notably, PJM's capacity auction in July 2024 resulted in record-high capacity prices of \$$270$/MW-day \cite{PJM1}. These charges are primarily tied to procuring sufficient generation capacity to meet peak demand. By managing usage during peak periods, commercial and industrial customers can potentially reduce their NCP and CP demand charges and thereby realize significant on-bill savings by optimally planning distributed energy resources such as on-site generation, battery energy storage systems (BESS), and controllable loads.

The present work focuses on the reduction of CP and NCP demand charges for commercial and industrial customers with microgrids. We formulate and solve a typical \emph{demand-response} problem as a stochastic optimization program that determines the next-day BESS schedule, minimizing the expected costs, and providing maximum on-bill savings. We model a microgrid that comprises a building responsible for the electricity demand, a local renewable generator in the form of solar panels, and a battery energy storage system (BESS). We use historical data from the demand and the solar energy production; we complement them with publicly available data of local weather as well as loads (and their forecasts whenever available) from the zone and the region of the microgrid. Using technology developed in \cite{carmona2024coincidentpeakpredictioncapacity} by some of the present authors, we produce 1) photovoltaic production forecasts, 2) Monte Carlo scenarios for the utility network load, 3) Monte Carlo scenarios for the local load, and 4) Monte Carlo estimates of the probabilities that CPs and NCPs will occur the next day (and at what times of the day they will occur). This allows us to fully specify the stochastic optimization problem. To provide numerical evidence of the performance of our algorithm, we work with data from an existing microgrid, which we \emph{anonymized} to protect the identity of the entity, which agreed to share information on the characteristics, the operation and the costs of the microgrid.  

To the best of our knowledge, this is one of the very first works to propose a unifying approach to address both CP and NCP demand charges, while taking into account the contribution of both a local PV generator and a BESS. Indeed, previous literature separately studied the reduction of NCP demand charges and that of CP charges. For instance, many approaches \cite{darghouth_demand_2020,ordonez_demand_2023} develop a methodology that evaluates the impact of photovoltaics and batteries on commercial and industrial load profiles and on the NCP demand charge savings depending on the charge program design. Other work \cite{Wu_Ma_Fu_Hou_Rehm_Lu_2022} focuses on the reduction of CP charges from the dispatch of a battery only, without any local generation. Similarly, approximate dynamic programming is used to address the mitigation of CP charges for data centers  \cite{Zhang19}. Management of the demand-side can be optimized through the scheduling of electric vehicle (EV) charging stations for demand charge reduction  \cite{Zhang18}, instead of looking at supply-side flexibility. We also notice methods to address the problem of optimal BESS sizing and scheduling in the Ontario, Canada jurisdiction \cite{kadri_energy_2020}: although the economic impact on demand charges is evaluated, it is not part of the optimization algorithm. There is also a rich body of work that focuses on BESS scheduling and optimization for participation and trading in multiple energy and ancillary services markets \cite{BESS1,SSN23}. However, these works are mostly attempting to optimize the bidding of storage facilities and they often ignore the management of peak demand charges that can provide significant cost savings to customers.

The rest of the paper is organized as follows. \autoref{sec:2} describes the setting of the coincident peak and demand charges programs and the data associated with the microgrid. It also describes the algorithms that are used to generate scenarios, predict the photovoltaic energy production, and estimate the CP and NCP peak probabilities. \autoref{sec:3} formulates the optimization problem of the battery. \autoref{sec:4} conducts a comparative analysis of the performance of the proposed methodology with a benchmark. 
 
\section{Model} \label{sec:2}

\subsection{Framework} 

\begin{figure}[htb!]
\centering
\begin{tikzpicture}[node distance=2cm, every node/.style={draw, fill=none, text centered, scale = 0.75}, font=\small]

    \tikzstyle{utility} = [rectangle, rounded corners, minimum width=2cm, minimum height=1cm, text centered, draw=black, fill=gray!40]
    \tikzstyle{controller} = [rectangle, rounded corners, minimum width=2cm, minimum height=1cm, text centered, draw=black, fill=teal!40]
    \tikzstyle{microgrid} = [rectangle, rounded corners, minimum width=8cm, minimum height=7cm, text centered, draw=teal,  dashed, fill=none]
    
    \tikzstyle{component} = [rectangle, rounded corners, minimum width=2cm, minimum height=0.8cm, text centered, draw=none, fill=none]
    
    \tikzstyle{load} = [rectangle, rounded corners, minimum width=1.5cm, minimum height=0.8cm, text centered, draw=black, fill=red!60]
    \tikzstyle{battery} = [rectangle, rounded corners, minimum width=1.5cm, minimum height=0.8cm, text centered, draw=black, fill=blue!40]
     \tikzstyle{solar} = [rectangle, rounded corners, minimum width=1.5cm, minimum height=0.8cm, text centered, draw=black, fill=orange!60]
    \tikzstyle{connection} = [draw=black, thick, -]
    \tikzstyle{arrow} = [->, thick]
    \tikzstyle{darrow} = [<->, thick]

    \node (utility) [utility] {Utility Grid};
    \node (controller) [controller, right of=utility, xshift=4cm] {Controller};
    \node (solar) [solar, above of=controller, yshift=1cm] {Solar Panels};
    \node (battery) [battery, below of=controller, yshift=-1cm] {Battery Storage};
    \node (load) [load, right of=controller, xshift=2cm] {Load};
    \node (microgrid) [microgrid, right of=load, xshift=-4.9cm] {};
    \node (component) [component, color=teal, right of=solar, xshift=2cm] {Microgrid};

    \draw[connection] (utility) -- (controller);
    \draw[connection] (controller) -- (solar);
    \draw[connection] (controller) -- (battery);
    \draw[connection] (controller) -- (load);

    \draw[darrow] (utility) -- (controller) node[midway, above, draw=none] {Power Flow};
    \draw[arrow] (solar) -- (controller) node[midway, left, draw=none] {Generated Power};
    \draw[darrow] (battery) -- (controller) node[midway, left, draw=none] {Stored Power};
    \draw[arrow] (controller) -- (load) node[midway, above, draw=none] {Supply};

\end{tikzpicture}
\caption{Diagram of the microgrid}
\label{fig:microgrid}
\end{figure}
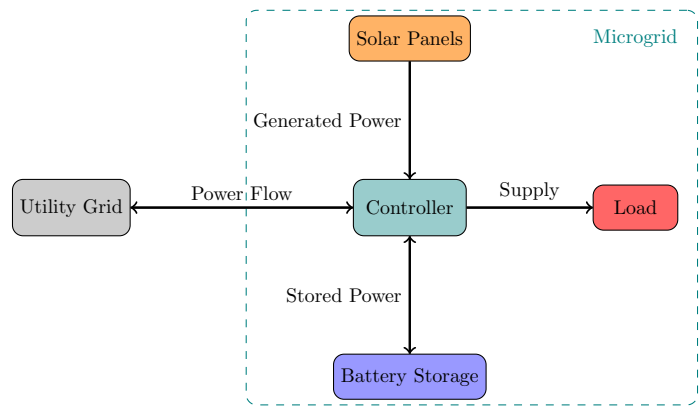

We consider a microgrid that is connected to the main grid through a utility provider and that comprises a load component (e.g. a building), a local renewable generator (e.g. solar panels), and a battery energy storage system (BESS) as depicted in \autoref{fig:microgrid}. 

The microgrid's load is assumed to be non-controllable: it has to be satisfied at all times and  cannot be modified by the controller (e.g. through reducing non-essential load generated by HVAC equipment or lighting). We assume that the local renewable generator consists of solar panels that produce electricity that can be used towards the load, to charge the BESS or to be sold back to the grid. 

\subsubsection{Pricing Structure, Capacity-Transmission Charges, and Demand Charges} 

There are two main compensation schedules for solar energy: net-metering will compensate at a flat rate that is the same when purchasing or selling energy, while net-billing will include a spread so that the price of selling energy to the utility is lower than the price to buy it. We place ourselves in the net-metering case and for simplicity, assume that it is a unique flat rate. 

\begin{remark}
    \begin{itemize}
        \item Although we assume that the rate is flat, our framework also accommodates a Time-Of-Use (TOU) schedule, which describes the adjustments to the price of electricity based on the time of day (e.g. on-peak, mid-peak, off-peak) and the season. 
        \item The net-billing case can also be addressed in our framework provided that some forecasts for the electricity prices are available. This will require adding an additional non-linear term in the objective function to account for the spread between selling and purchasing electricity.  
    \end{itemize}
\end{remark}

We assume that the microgrid is in the jurisdiction of the utility company PSE\&G and the Independent System Operator (ISO) PJM. It is therefore subject to several coincident-peak-based charges: 
\begin{itemize}
    \item Capacity charges through PJM's 5CP program, which records the load during the five highest ISO-wide hourly demand peak on different business days from the summer season, 
    \item Transmission charges through PSE\&G's 1CP program, which records the load during the highest utility-wide hourly demand peak on business days during the summer period. 
\end{itemize}

In practice, the transmission charges vary greatly from one utility to another, but for PSE\&G, they are significantly higher and therefore exceed the capacity charges by a non-unit factor for the same peak load contribution. To simplify the presentation of the algorithm, we focus on the PSE\&G single coincident peak program. The adaptation for PJM's 5CP program will be described where appropriate.

One additional thing to note is that the microgrid's own peak does not necessarily coincide with any of the ISO or the utility company's peak, that is why, we refer to the idiosyncratic peak as the non-coincident peak (NCP). The microgrid is subject to demand charges that are based on this non-coincident peak, which is the maximum amount of power drawn over a given time interval during a pre-determined timeframe. In PSE\&G's jurisdiction, the demand charge is computed on a monthly basis based on the microgrid's own hourly peak, regardless of what happens on the utility network.

\subsection{Input of the model}

\subsubsection{ISO and utility load}
The coincident peak prediction relies on the electric loads of two zones of interest $z_1$ and $z_2$ with the characteristic that zone $z_2$ is a subset of zone $z_1$. We assume that the CP program is ran at the level of zone $z_2$ and that we have the following information at an hourly resolution: 
\begin{itemize}
    \item the actual load historical data of the region $z_1$,
    \item the day-ahead forecast data  of the region $z_1$,
    \item the actual load historical data of the zone $z_2$.
\end{itemize}

The ISO PJM publishes these datasets for $z_1 = \texttt{MA}$ and  $z_2 = \texttt{PS}$ where the region Mid-Atlantic, denoted \texttt{MA}, regroups utilities in the corresponding geographic area, including the local utility company PSE\&G, denoted \texttt{PS}. To simplify the presentation, we consider a single forecast with a forecasting horizon that covers all 24 hours of the day the predictions are generated for. 

\subsubsection{Microgrid's load, local generation and battery management}

The use case of our study is an anonymized site that is served by PSE\&G and falls under the jurisdiction of PJM. The historical data available for this microgrid covers the period between 2022-05 and 2023-12 with a 15-minute frequency, which we resample to an hourly frequency. It consists in the building historical load, the historical photovoltaic (PV) power production and the historical battery schedule. 

\subsubsection{Meteorological data}

Our battery management optimization algorithm requires forecasts for the hourly photovoltaic (PV) power production over the next day. Since such forecasts are not publicly available, we produce them with a regression model fitted to historical meteorological data from the ERA5 dataset \cite{era5} and historical PV productions at the location of interest. The meteorological data include hourly forecasts of temperature, cloud cover, and solar radiation, as generated the day before.

\subsubsection{Battery specifications}

For the sake of completeness, we provide the technical specifications of the Battery Energy Storage System in \autoref{tab:tech_specs_bess}. 

\begin{table}[htb!]
    \centering
    \begin{tabular}{|c|c|c|}
    \hline 
       Notation  &  Description & Value\\ \hline 
       $C$  & Maximum capacity & $1000$ MWh\\
       $\Pi^+$ & Maximum discharging output & $500$ MW \\ 
       $\Pi^-$ & Maximum charging output & $500$ MW\\ 
       $\eta$ & Round-trip efficiency & $0.8$ \\
       $\eta^+$ & Discharging efficiency & $\sqrt{0.8}$\\ 
       $\eta^-$ & Charging efficiency & $\sqrt{0.8}$\\ 
       $Q_h$ & Available capacity at hour $h$ & N/A \\
       \hline 
    \end{tabular}
    \caption{Technical specifications of the BESS}
    \label{tab:tech_specs_bess}
\end{table}

\subsubsection{Notations} We will denote by $\Iintv{i,j}$ the set of integers comprised between $i$ and $j$ such that $i,j \in \mathbb{N}$ and $i\leq j$. 

\subsection{Building blocks} \label{ssec:building_blocks}

\subsubsection{CP Prediction}

The prediction of coincident peak events relies on the Monte-Carlo scenario generation engine, \texttt{PLProb}, developed in \cite{carmona2024coincidentpeakpredictioncapacity}, publicly available for academic research purposes. It can generate scenarios for the electric load of the jurisdiction running the coincident peak programs in two cases, depending on the availability of a direct forecast or not. The algorithm also produces estimates for the prediction of CP-day and CP-hour probabilities, corresponding respectively to the probability that the next day sees the running coincident peak and the probability at each hour that this CP event may happen. 

The version used to predict PSE\&G's single CP is the conditional scenario version that is able to create scenarios for PSE\&G ($\texttt{PS}$) based on the Mid-Atlantic load data ($\texttt{MA}$). We provide a stylized description of the algorithm below and refer the interested reader to \cite{carmona2024coincidentpeakpredictioncapacity} for the technical details. 

\begin{center}
    \captionof{algorithm}{Conditional Utility Load Scenario Generation}
    \label{algo:cond_scen}

\begin{enumerate}
    \item Using historical data, build a joint model for the actual loads of $\texttt{MA}$ and their forecasts, and use this model to generate Monte Carlo scenarios of the $\texttt{MA}$ actual loads from a set of their current forecasts;
    \item Using historical data again, build a joint model for the actual loads of $\texttt{MA}$ and $\texttt{PS}$;
    \item Use this model and the scenarios of the $\texttt{MA}$ actual loads
    constructed in the first step to generate $N_s=1000$ scenarios for the PSE\&G loads by conditional Monte Carlo simulation;
    \item Use these scenarios to compute the probability that tomorrow is going to be a new running maximum hourly load day, call it $p^{\texttt{PS}}_{d,nrm}$;
    \item Use the same scenarios to compute for each hour $h\in \Iintv{0,23}$ the probability that the PSEG load on hour $h$ will be the highest hourly load of the day, and call this probability $p^{\texttt{PS}}_h$.
\end{enumerate}
\end{center}

\autoref{fig:proba_histo_example} shows the histogram of the probability of occurrence of the highest daily load at each hour for two different days. 

\begin{remark}
    Predicting PJM's five CPs relies on the standard scenario engine that requires the actual load and the forecast data published by PJM itself. Similarly, we can generate the Monte Carlo scenarios and estimate the probabilities that the next day will see one of the five highest running maxima of PJM's load and the probability of each hour to see this event. We refer to \cite{carmona2024coincidentpeakpredictioncapacity} for details. 
\end{remark}

\subsubsection{Photovoltaic Power Production Forecast}

\begin{figure}[htb!]
    \centering
     \includegraphics[width=0.6\linewidth]{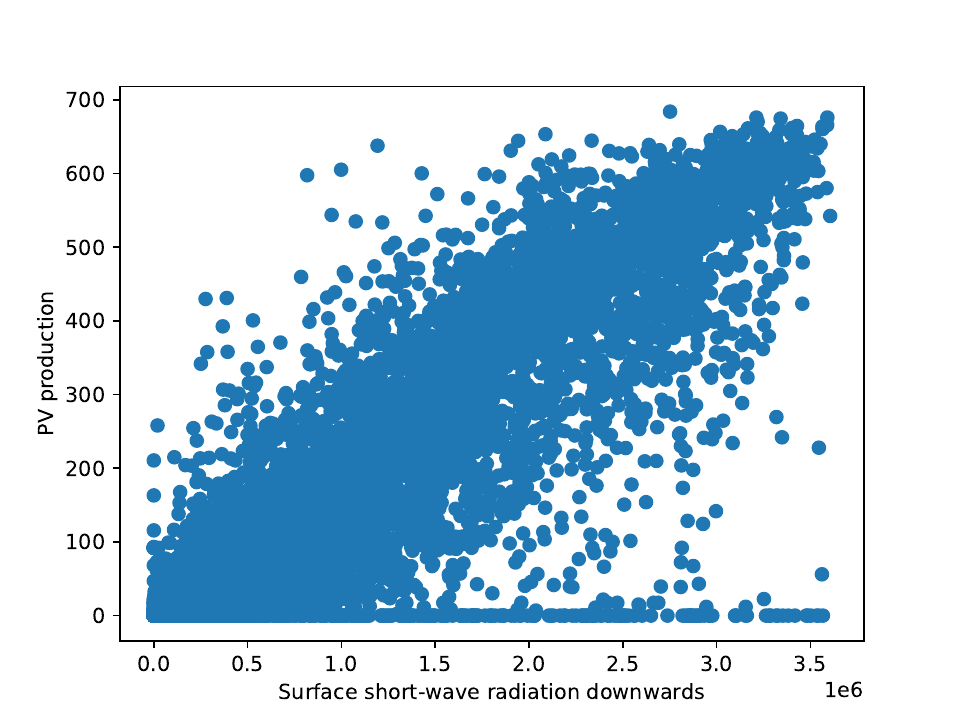} 
    \caption{Scatter plot of surface short-wave radiation downwards versus PV production.}
    \label{fig:ssrd_pv}
\end{figure}

The primary objective of this study is to present an optimization framework for battery management. While accurate PV forecasts are crucial, developing a highly predictive forecast model is beyond the scope of this study. Therefore, we opt for a simple linear regression model. Specifically, we use a non-intercept linear regression model with one predictor variable: surface short-wave radiation downwards. \autoref{fig:ssrd_pv} shows a scatter plot of surface short-wave radiation downwards versus PV production. The horizontal ``line'' near zero PV output is possibly due to PV maintenance. To account for variations in the relationship between the response and predictor variables at different hours, we separate the training data by hour. In summary, for a given day and hour, we fit a linear regression model using all available historical forecast data for that hour (prior to the day in question) and the corresponding historical PV production. Once the linear coefficient is determined, a PV prediction is made by multiplying the forecasted surface short-wave radiation downwards by the coefficient. The mean absolute error (MAE) and mean squared error (MSE) of our models for the period from 2023-01-01 to 2023-12-31 are respectively 35.2 and 5,264.9.

\subsubsection{Local Load Scenario} 
We denote by $\texttt{MG}$ the microgrid building load and by  $\texttt{MG-net}$ the microgrid net load (net of photovoltaic production and battery load). We will generate Monte Carlo scenarios for the microgrid building load and the microgrid net load conditionally on the scenarios  previously generated  for $\texttt{PS}$ by Algorithm \ref{algo:cond_scen}. 

\begin{center}
    \captionof{algorithm}{Conditional Local Load Scenario Generation}
    \label{algo:cond_local_load_scen}
    \begin{enumerate}
        \item Use historical data to build a joint model for the actual loads of $\texttt{PS}$ and $\texttt{MG}$;
        \item Use this model and the Monte Carlo scenarios of the actual $\texttt{PS}$ load generated using Algorithm \ref{algo:cond_scen}, to produce $N_s=1000$ scenarios for the local $\texttt{MG}$ load for the next day. Call them $L^{\texttt{MG}}_1,L^{\texttt{MG}}_2,\cdots,L^{\texttt{MG}}_{N_s}$ where each $L^{\texttt{MG}}_s$ is a vector of the 24 hourly load for scenario $s \in \Iintv{1, N_s}$.
    \end{enumerate}
\end{center}

\begin{figure}[htb!]
    \centering
    \includegraphics[width=0.45\linewidth]{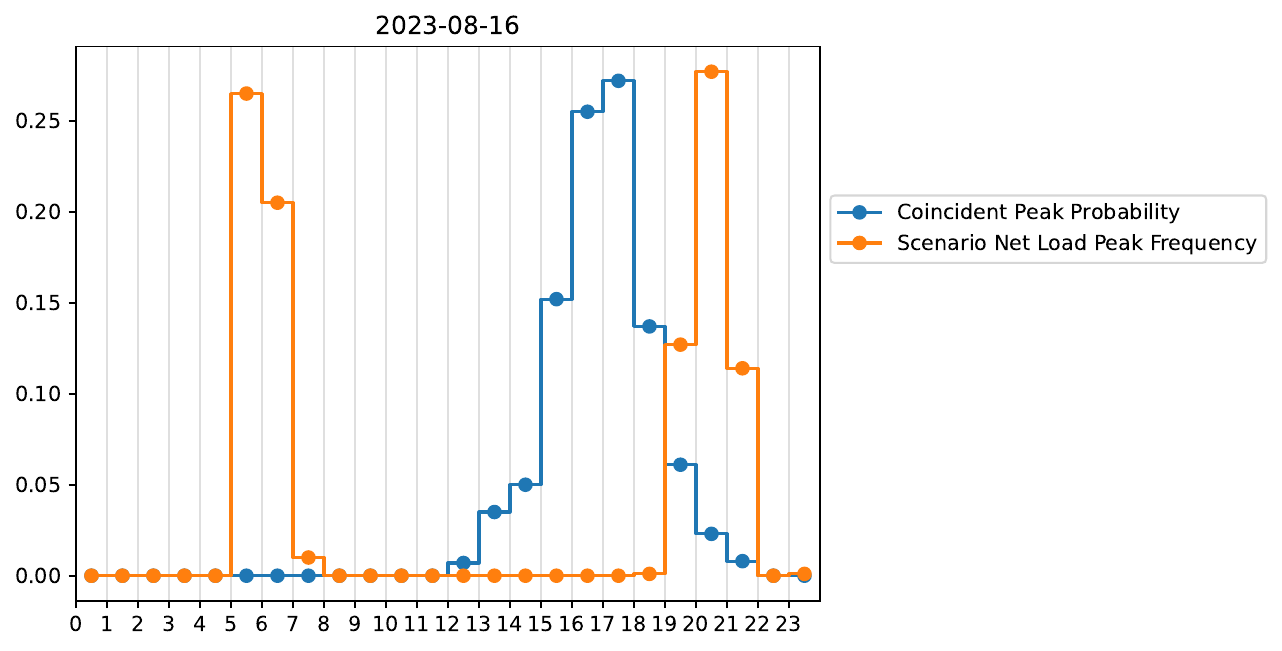} 
     \includegraphics[width=0.45\linewidth]{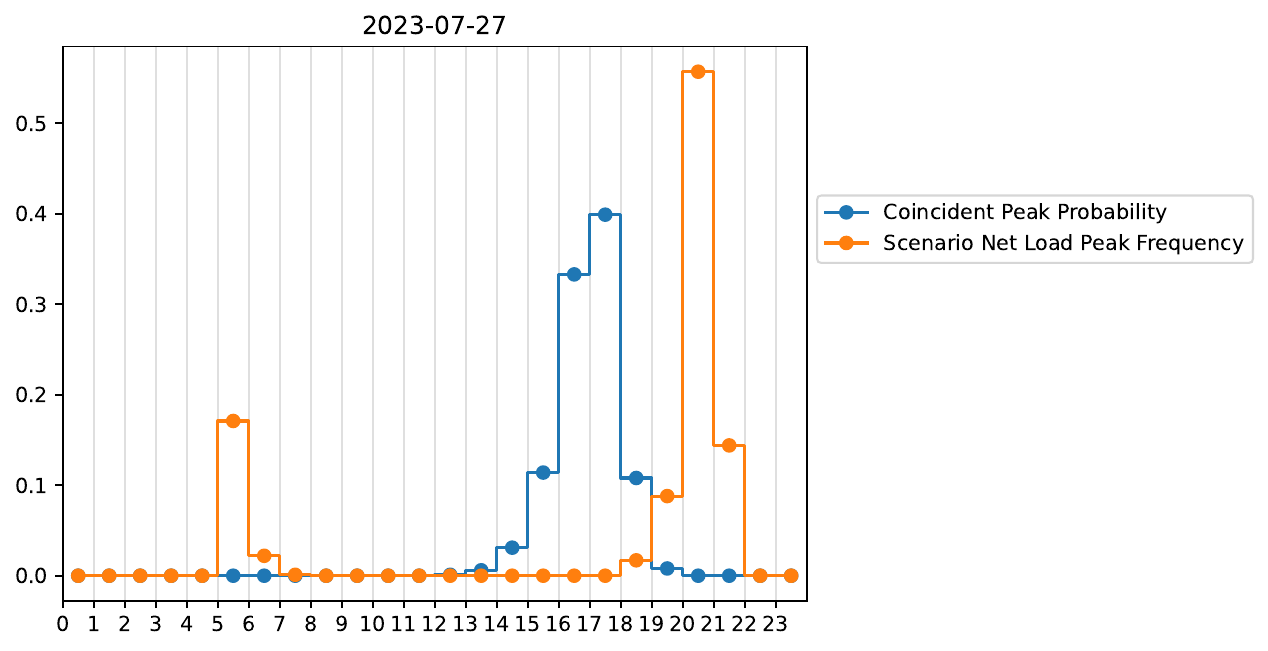} 
    \caption{Hourly probability histogram of the running CP daily maximum (in blue) and the running NCP daily maximum (in orange) for \texttt{2023-08-16} and \texttt{2023-07-27}: the NCP histogram shows the probability estimates prior to any optimization and therefore excludes any additional load originating from the battery.}
    \label{fig:proba_histo_example}
\end{figure}

The scenarios for the microgrid building load will be leveraged during the formulation of the optimization problem. On the other hand, the non-coincident peak (NCP) program penalizes the hourly peak load of the microgrid for the current month. To account for this program, we apply Algorithm \ref{algo:cond_local_load_scen} to generate $N_s=1000$ scenarios of the microgrid's net load $\{\mathbf{L}^{\texttt{MG-net}}_{1},\cdots,\mathbf{L}^{\texttt{MG-net}}_{N_s}\}$, where each component $\mathbf{L}^{\texttt{MG-net}}_{s}$ is a vector of $24$ hourly loads. From the history of the microgrid load in the current month before day $d$, and these $N_s$ scenarios, we derive estimates, $p^{\texttt{MG}}_{d,nrm}$, of the probability that day $d$ sees a new running maximum, i.e. a candidate for being the NCP-day of the month. From the scenarios $L^{\texttt{MG}}$ and the PV production forecasts, we plot the frequency of occurrence of maximum net load excluding any contribution from the battery in \autoref{fig:proba_histo_example}. 

\subsubsection{Benchmark}
To protect the anonymity of the microgrid whose historical data we use, we compare the results of our optimization method to those of a baseline battery management heuristic based on rational principles. Like our algorithm, this benchmark is based on historical data: the strategy consists in charging and discharging the battery on a schedule to minimize the demand charges and discharging on a fixed time window to catch the CP only when an alert is received.
\begin{itemize}
    \item There are two main peaks for the power exchanged with the grid: one in the morning and one in the evening. 
    \item The CP-day alert will be based on a signal from the MIDATL forecast published by PJM. 
    For the test year 2023, the hyperparameters for the signal will be tuned based on the backtest performance from 2014 up until 2022 (included). 
    \item The PSE\&G 1CP has been occurring between 4PM and 6PM between 2015 to 2022. 
\end{itemize}

\begin{figure}[htb!]
    \centering
    \begin{subfigure}[t]{0.45\textwidth}
        \includegraphics[width=\textwidth]{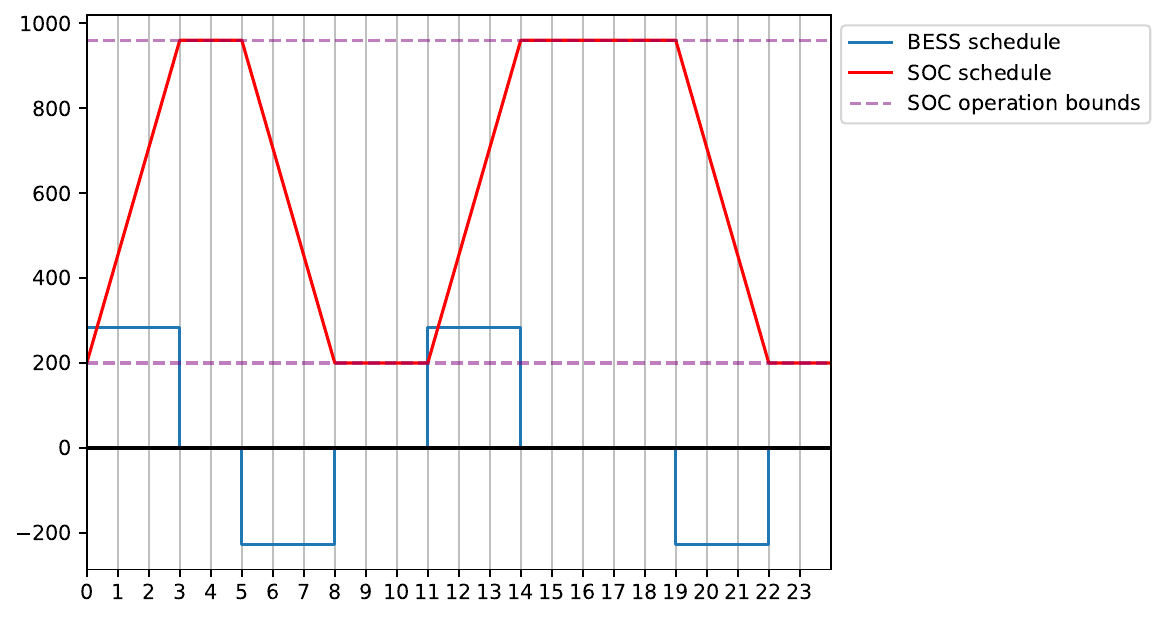}
        \caption{With no CP alert}
    \end{subfigure}
     \begin{subfigure}[t]{0.45\textwidth}
        \includegraphics[width=\textwidth]{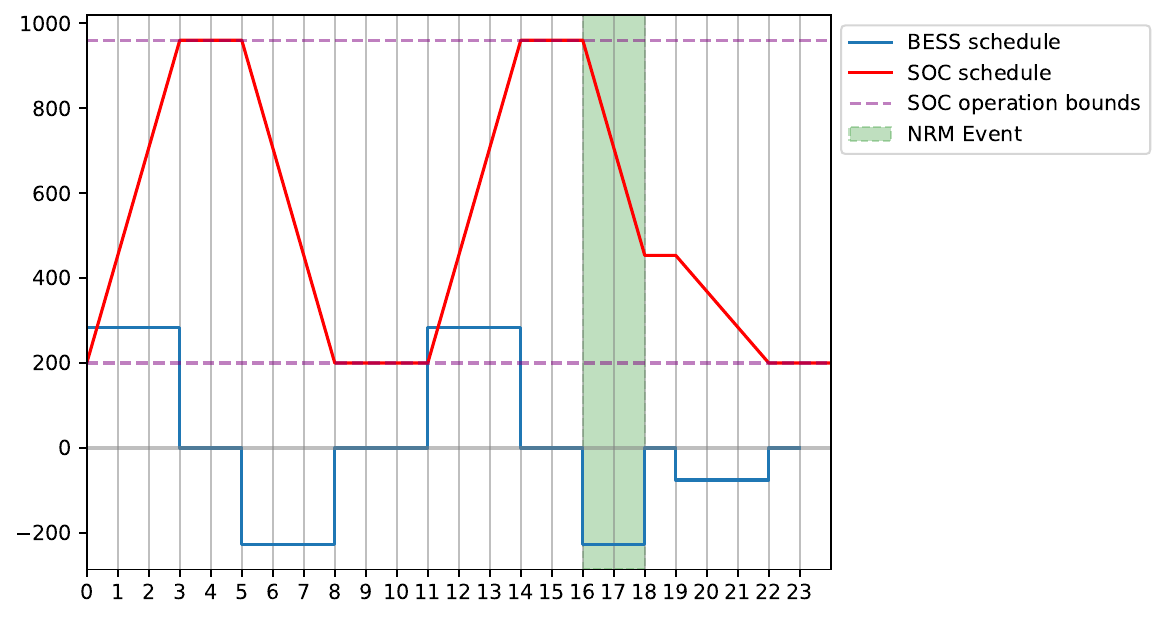}
        \caption{With a CP alert}
    \end{subfigure}
    \caption{Benchmark BESS Schedule}
\end{figure}

Using the data of the building load net of PV production from Summer 2022, a rule-based schedule for 2023 can be determined as the following:

\begin{itemize}
    \item Charge at night between 12AM-3AM at a constant rate when the load is usually low; 
    \item Discharge to shave Peak 1 at a constant rate between 5AM-8AM 
    \item Charge at a constant rate between 11PM-2PM,  when the load is usually low;
    \item If there is a CP alert, discharge at the discharge rate between 4PM-6PM; 
    \item Discharge what is still available to shave Peak 2 that usually occurs 7PM-10PM. 
\end{itemize}

\begin{remark}
    This benchmark relies on a prediction of the CP events, solely based on the publicly available data. In practice, commercial and industrial customers pay for the service of external consulting companies to get such predictions. This often expensive cost will not be included in the economics of the comparison. 
\end{remark}

\section{Battery Optimization Problem} \label{sec:3}

\subsection{Control variables, constraints, and admissible solutions}

Let $Q_h$ be the available capacity of the battery at hour $h$, i.e. $Q_h = C \times SOC_h$ where $SOC_h$ stands for State of Charge that takes values between 0\% and 100\%. Let $\pi^{B}_h$ be the power used to charge or discharge the battery at hour $h$ with $\pi^{B}_h \in [-1,1]$. The convention is that $\pi^{B}_h \geq 0$ is for discharging the BESS and $\pi^{B}_h \leq 0 $ is for charging it. We introduce $\pi^+_h, \,\pi^-_h \in[0,1] $ so that $\pi^{B}_h = \pi^+_h - \pi^-_h$ for every hour $h$ and we denote by $\bpi^B_h = (\pi^-_h, \pi^+_h)^{\dagger}$ and $\bpi^B=(\pi^-_1, \pi^+_1,\cdots,\pi^-_{24}, \pi^+_{24})$ the controls. 

By taking efficiency into account, the physical constraints of the battery incorporate: 
\begin{itemize}
    \item An initial SOC of $SOC_0 = 20 \%$;
    \item Maximum rate of charge/discharge constraints;
    \begin{align*}
        & 0 \leq \pi^-_h C \leq  (1-b_h) \Pi^- & \text{(charging)} \\
        & 0 \leq \pi^+_h C \leq  b_h \Pi^+, & \text{(discharging)} 
    \end{align*}
     where the binary variables $b_h$ ensure that charging or discharging occurs one at a time;
    \item The dynamical equation for the SOC as $ SOC_h=SOC_{h-1}+ \pi^-_h\eta^- - \pi^+_h/\eta^+ $
    where $\eta^-$ and $\eta^+$ are the charging and discharging efficiencies. 
    \item Preservation constraints: $\underline{SOC} \leq SOC_h \leq \overline{SOC}, \forall h \in \{1, \dots, 24\}$ where $\underline{SOC} = 20 \%$ and $\overline{SOC} = 96 \%$. 
\end{itemize}

We say that $\bpi^B=(\pi^{-}_1,\pi^{+}_1,\cdots,\pi^{-}_{24},\pi^{+}_{24}, )$ is an \textbf{admissible} strategy profile, in notation $\bpi\in\mathcal{A}$, if for every $h\in \Iintv{0, 23}$, $\pi^{+}_h, \, \pi^{-}_h \in[0,1]$ and if the above physical constraints are satisfied. 

\subsection{Expected P\&L}

Assume momentarily that $\pi\in\mathcal{A}$ is fixed. 

\subsubsection{Daily Electricity Consumption Cost}

For each scenario $i\in\{1,\cdots,N_s\}$ we compute the electricity consumption cost (ECC) for scenario $i$ as
$$
\operatorname{ECC}_i= \sum_{h=0}^{23}P_h \times \Big(L^{\texttt{MG}}_{i,h} - PV_h - C (\pi^+_h - \pi^-_h) \Big)
$$
and the expected total consumption cost as
$$
\operatorname{ECC} = \frac1{N_s}\sum_{i=1}^{N_s} \operatorname{ECC}_i
$$

\subsubsection{Daily Running Coincident Peak P\&L}
We consider the expected reward from the CP program given by
$$
\operatorname{ECP} =  \lambda_{CP} \cdot p^{\texttt{PS}}_{d,nrm}\sum_{h=0}^{23}p^{\texttt{PS}}_h \Big( \bar L^{\texttt{MG}}_{h} - PV_h - C(\pi^+_h - \pi^-_h)\Big)
$$ 
where $\bar L^{\texttt{MG}}_{h}=\frac{1}{N_s}\sum_{i=1}^{N_s}L_{i,h}$ is the expected local load (computed over the $N_s$ scenarios) at hour $h$ and $\lambda_{CP}$ is a coefficient associated with the capacity and transmission charges. In practice, $\lambda_{CP}$ is the product of a zonal factor (specific to the microgrid location), 30 (for 30 days) and a NITS rate in $ \$/\text{Day}/MW$ (this rate is established from planned infrastructure construction work). 

\begin{remark}
    Including the capacity component reward requires the addition of a term of the following form: 
$$
\operatorname{ECP}^{\texttt{PJM}}  =  \lambda_{5\text{CP}} \cdot p^{\texttt{PJM}}_{d,nrm}\sum_{h=0}^{23}p^{\texttt{PJM}}_h \Big( \bar L^{\texttt{MG}}_{h} - PV_h - C(\pi^+_h - \pi^-_h)\Big)
$$ 
where the probabilities $p^{\texttt{PJM}}_{nrm}$ and $p^{\texttt{PJM}}_h$ are similarly defined but meant to capture one of the five running highest peaks of the PJM system load, and $\lambda_{5\text{CP}}$ will be the product of a zonal factor, 30 (for 30 days), a Capacity Price in  $\$/\text{Day}/MW$  and $\frac{1}{5}$ to account for the contribution of a single coincident-peak load to the Peak Load Contribution (which is the average on the five CPs). 
\end{remark}

\subsubsection{Daily Demand Component P\&L} 

The demand component P\&L, denoted $EDC$, corresponds to the reward earned by the microgrid through the non-coincident peak (NCP) program, which penalizes the hourly peak load of the microgrid for the current month. The associated reward takes the following form: 
$$
\operatorname{EDC} =  \lambda_{NCP} \cdot p^{\texttt{MG}}_{d,nrm}\frac{1}{N_s}\sum_{s=1}^{N_s} \Big( L^{\texttt{MG}}_{s,h^*} - PV_{h^*} - C(\pi^+_{h^*} - \pi^-_{h^*}) \Big)
$$
where $h^*$ is the hour of the maximum load on the grid for that day, i.e.
$$
h^*=\operatorname{argmax}_{h=1,\cdots,24}L^{\texttt{MG}}_{s,h} - PV_{h} - C(\pi^+_h - \pi^-_h).
$$

\subsubsection{Daily Battery Degradation Penalty}

For obvious reasons, we include a degradation penalty resulting for battery over-usage. We choose it to be a linear function of the Depth of Discharge (DOD) expressed through the SOC changes or the power entering and exiting the battery: 
$$ \operatorname{EDP}  =  \lambda_{\text{deg}} \sum_{h=0}^{23} \lvert \pi^-_h - \pi^+_h - \pi^{-}_{h-1}  + \pi^{+}_{h-1}\rvert $$
Finally, the optimal strategy profile $\pi^*$ is computed as 
$$
\bpi^*= \underset{\bpi \in \mathcal{A}}{\operatorname{argmin}} \;\; OBJ(\bpi).
$$
the objective function $OBJ$ being defined by:
$$
\operatorname{OBJ}(\bpi) =  \operatorname{ECC}+ \operatorname{ECP} + \operatorname{EDC} + \operatorname{EDP}.
$$
\begin{remark}
It is worth clarifying the roles played by the various Monte Carlo scenarios in the above computations and the final optimization algorithm. 
\begin{itemize}
    \item On any given day, scenarios for the PSE\&G hourly load are used to compute the probabilities $p^{\texttt{PS}}_{d, nrm}$ and $p^{\texttt{PS}}_h$ which are estimates of the likelihood that next day is going to be a CP day, and of the likelihoods that each of the next $24$ hours will be the hour the daily hourly maximum occurs. So should other estimates, say $p^{\texttt{MG}}_{d, nrm}$ and $p^{\texttt{MG}}_h$, be available, they could be fed into the above formulas, the optimization of the objective function $OBJ$ be performed with these alternative estimates, in which case we would not need to use scenarios for the PSE\&G hourly load.
    \item Similarly, the scenarios $L^{\texttt{MG}}_s$ for the hourly needs in electricity for the microgrid's building, only enter the computation of the objective function through their hourly average $\bar L^{\texttt{MG}}_h$. So should other estimates, say $\hat{\bar{L}}^{\texttt{MG}}_h$, be available for the hourly loads required by the microgrid's building the next day, it could be used into the above formulas, the optimization of the objective function $OBJ$ be performed with these alternative estimates, in which case we would not need to use scenarios for the microgrid hourly load either.
    \item So if forecasts for the probabilities and the hourly loads of the microgrid are available, the above optimization can be performed without having to rely on Monte Carlo scenarios. This is due to the linear nature of the formula giving the objective function $OBJ$. See below for a nonlinear objective function requiring the use of individual scenarios.
\end{itemize}    
\end{remark}

Linearizing the demand charge term and the absolute value in the degradation cost leads us to the following Mixed-Integer Linear Program (MILP) problem:

\begin{align*}
    \underset{\bpi}{\min} & \quad \overline{\operatorname{OBJ}}(\bpi) & \numberthis \label{eq:milp}\\ 
    \text{s.t.} & & \\
    & \operatorname{SOC}_h=\operatorname{SOC}_{h-1}+ \pi^-_h \eta^- - \pi^+_h/\eta^+, &h \in \Iintv{0,23}\\
    & 0 \leq \pi^-_h C \leq  (1-b_h) \Pi^-,&  h \in \Iintv{0,23}\\
    & 0 \leq \pi^+_h C \leq  b_h \Pi^+, & h \in \Iintv{0,23}\\
    & \underline{\operatorname{SOC}} \leq \operatorname{SOC}_h \leq \overline{\operatorname{SOC}},&  h \in \Iintv{0,23}\\ 
    & \pi^-_h  - \pi^+_h -\pi^{-}_{h-1} + \pi^{+}_{h-1}   = u_h - v_h, & h \in \Iintv{0,23}\\
    & s_i \geq L^{\texttt{MG}}_{i,h} - PV_h + C \pi^-_h - C \pi^+_h, &  i \in \Iintv{1,N_s},  h \in \Iintv{0,23}\\
    & 0 \leq \pi^+_h \leq 1,  & h \in \Iintv{0,23}\\
    & 0 \leq \pi^-_h \leq 1, &  h \in \Iintv{0,23}\\
    & b_h \in \{0, 1\}, & h \in \Iintv{0,23}\\
    & u_h \geq 0, &  h \in \Iintv{0,23}\\
    & v_h \geq 0,  &  h \in \Iintv{0,23}\\ 
    & s_i \geq 0, &  i \in \Iintv{1,N_s}
\end{align*}

where $ \overline{\operatorname{OBJ}}(\bpi) = \operatorname{ECC} + \operatorname{ECP} + \lambda_{deg} \sum_{h=0}^{23} \{ u_h + v_h\} + \lambda_{NCP}\cdot \hat{p}_{d, nrm} \frac{1}{N_s} \sum_{i=1}^{N_s} s_i$ 
with the convention that $\operatorname{SOC}_0 = \underline{\operatorname{SOC}}$, $\pi^+_0 = 0, \pi^-_0  = 0, \pi^{PV}_0 = 0$. 
MILP \ref{eq:milp} is solved with \texttt{Gurobi} \cite{gurobi}.


\section{Numerical results} \label{sec:4}

We illustrate the performance of our algorithm with results obtained for a set of days whose likelihoods of being a CP or NCP cover a wide spectrum of values. 

\subsection{Comparison of Daily Performance}

The CP and NCP probabilities of the days we chose to investigated are reported in \autoref{tab:cp_prob_date}. 
\begin{table}[htb!]
    \centering
    \resizebox{\textwidth}{!}{
    \begin{tabular}{|c|c|c|c|} 
    \hline 
        Date & CP-day probability & Maximum CP-hour probability & NCP proba \\ \hline 
        2023-06-30 & 0.504 & 0.274 & 0.155\\
        2023-07-27 & 0.998 & 0.394 & 0.155\\
        2023-09-07 & 0.962 & 0.362 & 0.768\\
        2023-09-21 & 0.00 & 0.287 & 0.199\\
        \hline 
    \end{tabular}
    }
    \caption{Probability of CP-events for selected dates where \texttt{2023-09-07} is the true final CP}
    \label{tab:cp_prob_date}
\end{table}

\begin{figure}[htb!]
    \centering
    \begin{subfigure}[b]{\textwidth}
    \includegraphics[width=0.45\linewidth]{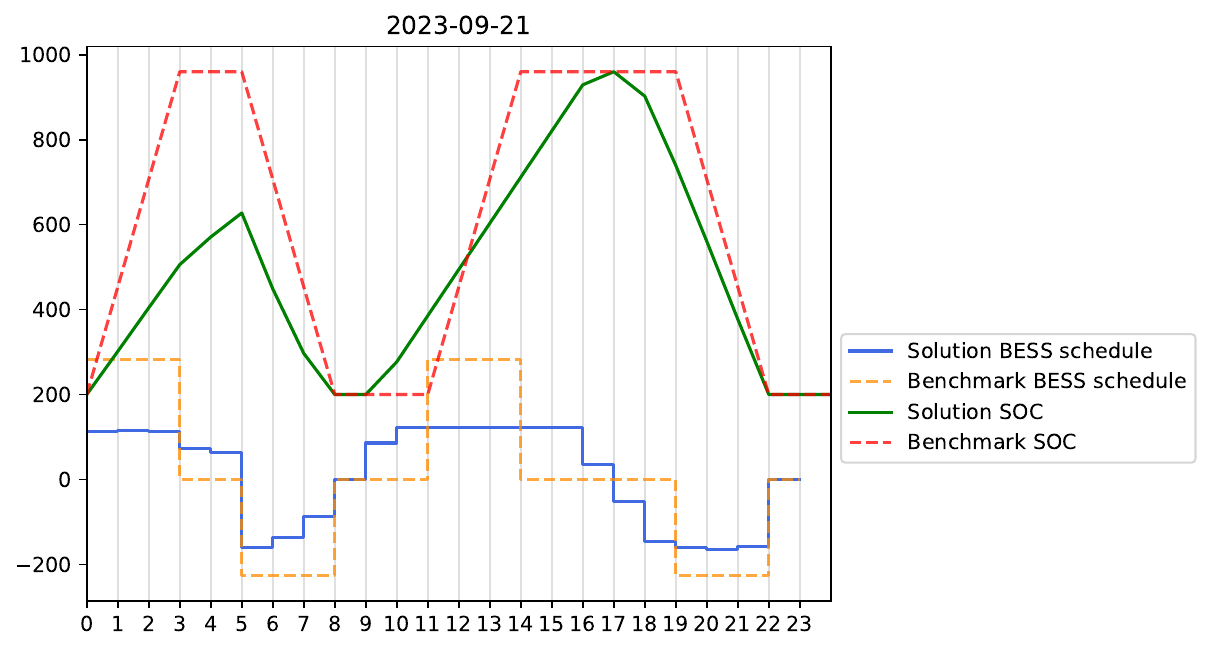} \hfill 
    \includegraphics[width=0.48\linewidth]{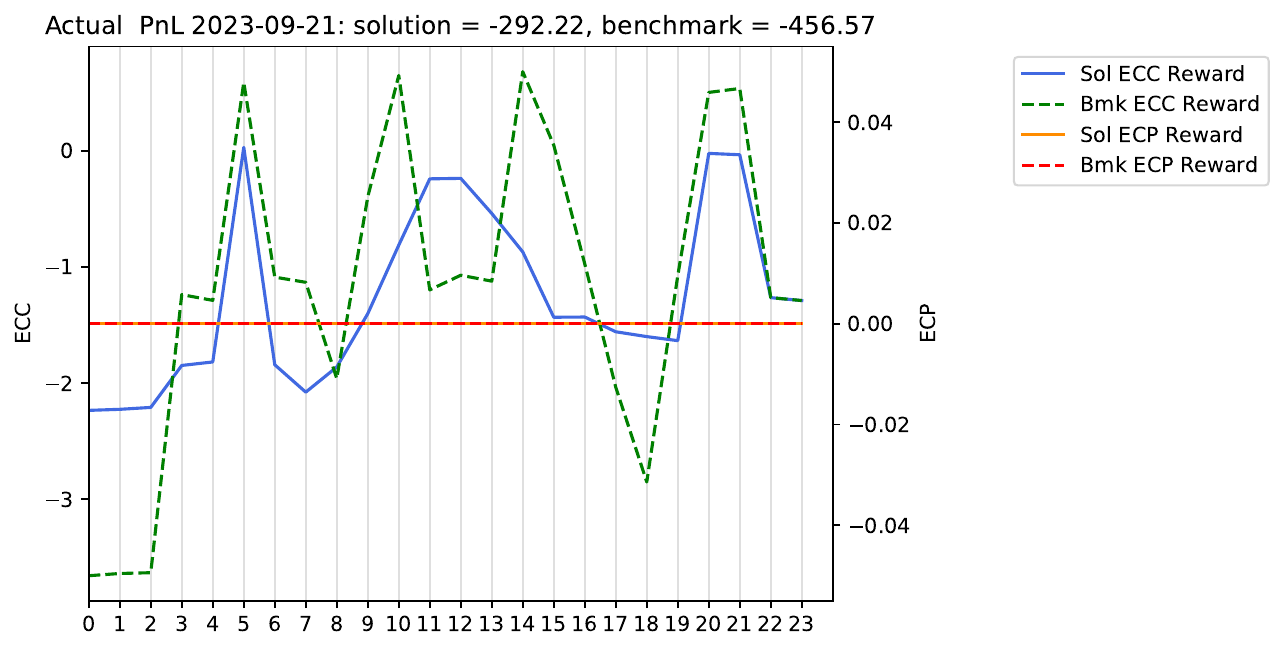}
    \caption{Performance on a day with zero CP-day probability, \texttt{2023-09-21}}
    \label{sfig:low_prob}
    \end{subfigure}
    
    \begin{subfigure}[b]{\textwidth}
   \includegraphics[width=0.45\linewidth]{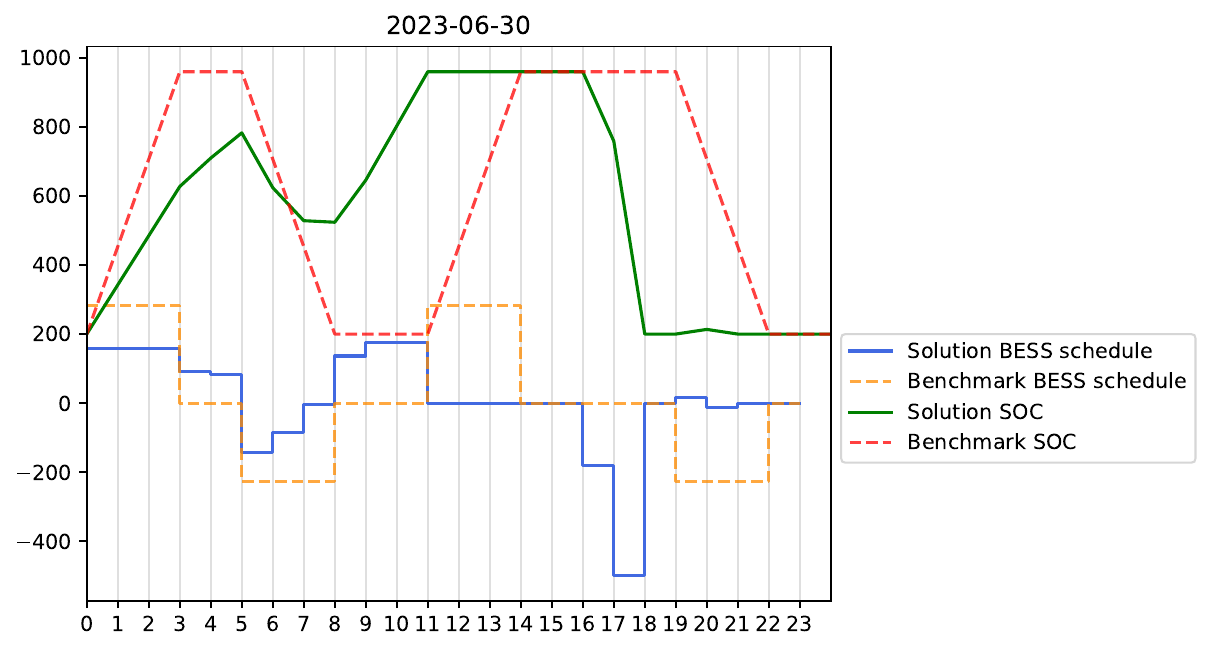} \hfill 
    \includegraphics[width=0.48\linewidth]{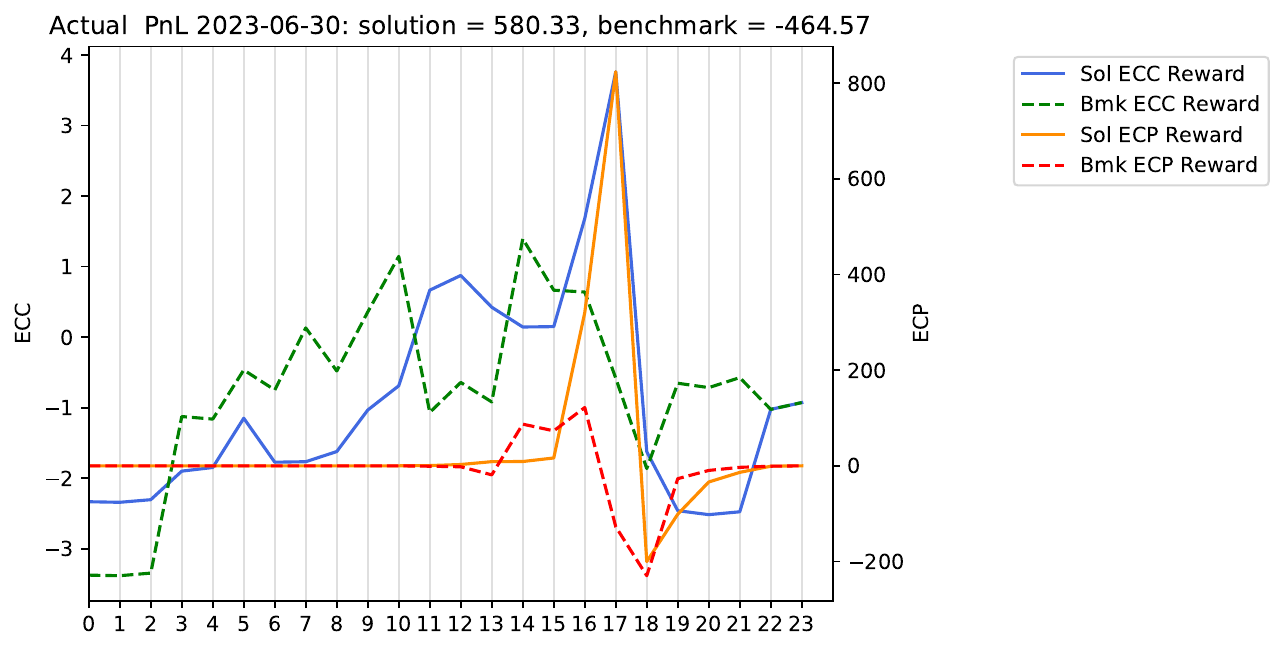}
    \caption{Performance on a day with medium CP-day probability, \texttt{2023-06-30}}
    \label{sfig:med_prob}
    \end{subfigure}
    
    \begin{subfigure}[b]{\textwidth}
    \includegraphics[width=0.45\linewidth]{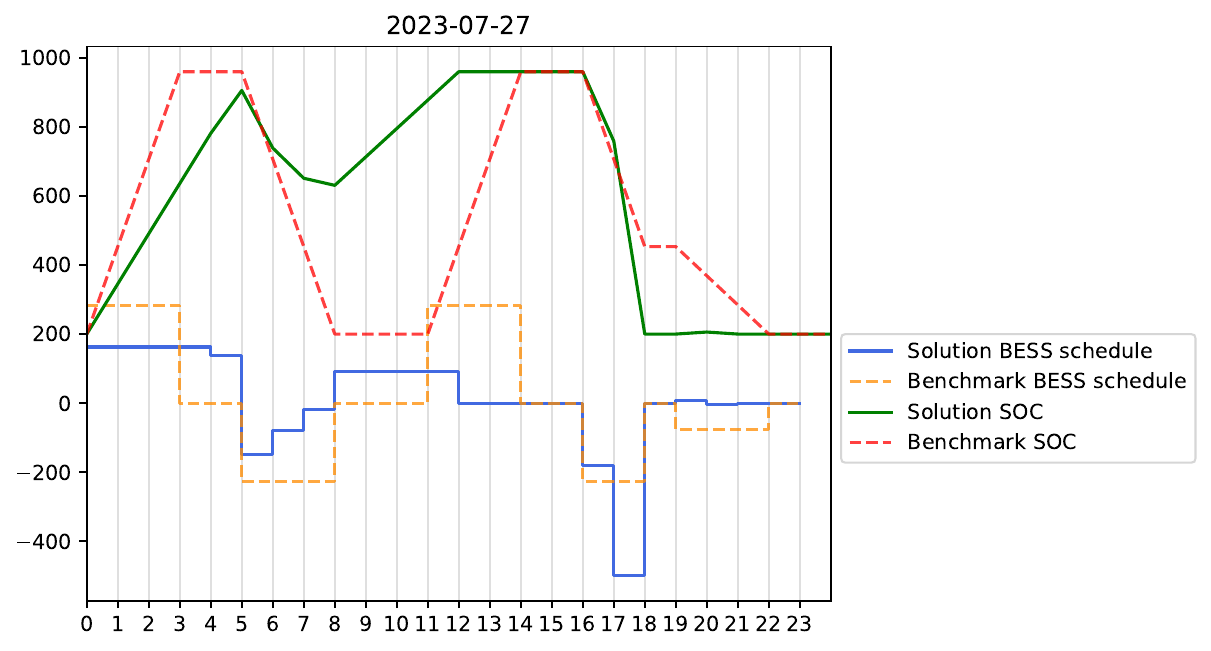} \hfill 
    \includegraphics[width=0.48\linewidth]{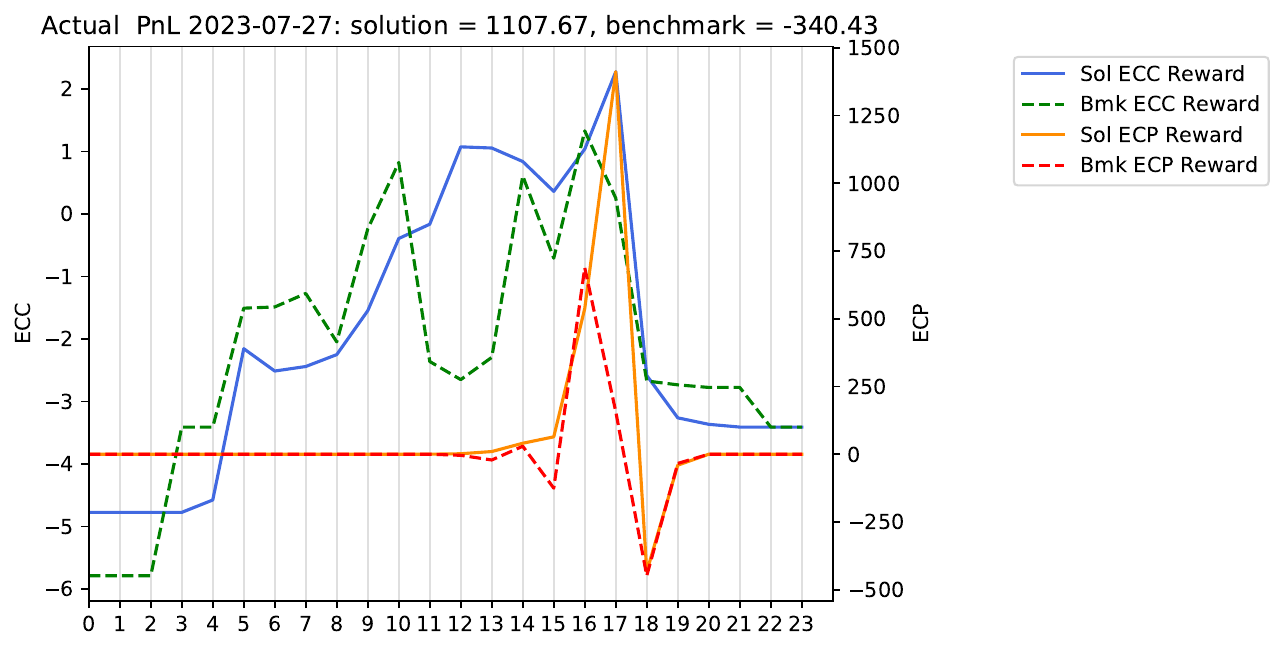}
    \caption{Performance on a day with high CP-day probability:, \texttt{2023-07-27}}
    \label{sfig:high_prob}
    \end{subfigure}

    \begin{subfigure}[b]{\textwidth}
    \includegraphics[width=0.45\linewidth]{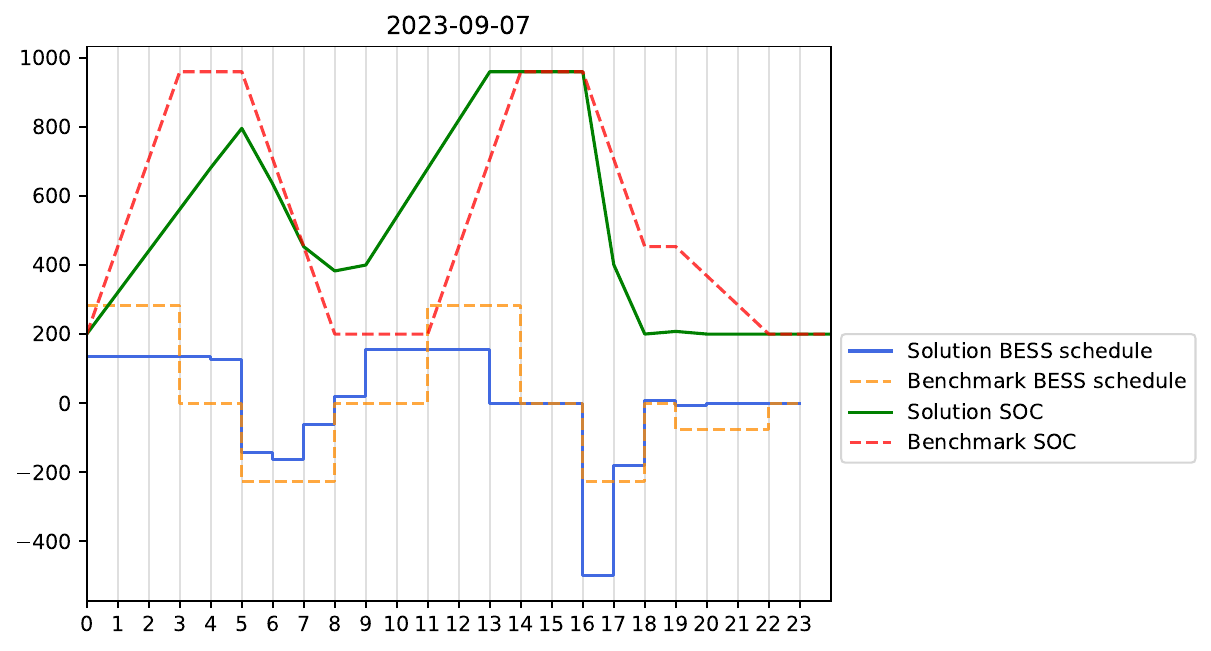} \hfill 
    \includegraphics[width=0.48\linewidth]{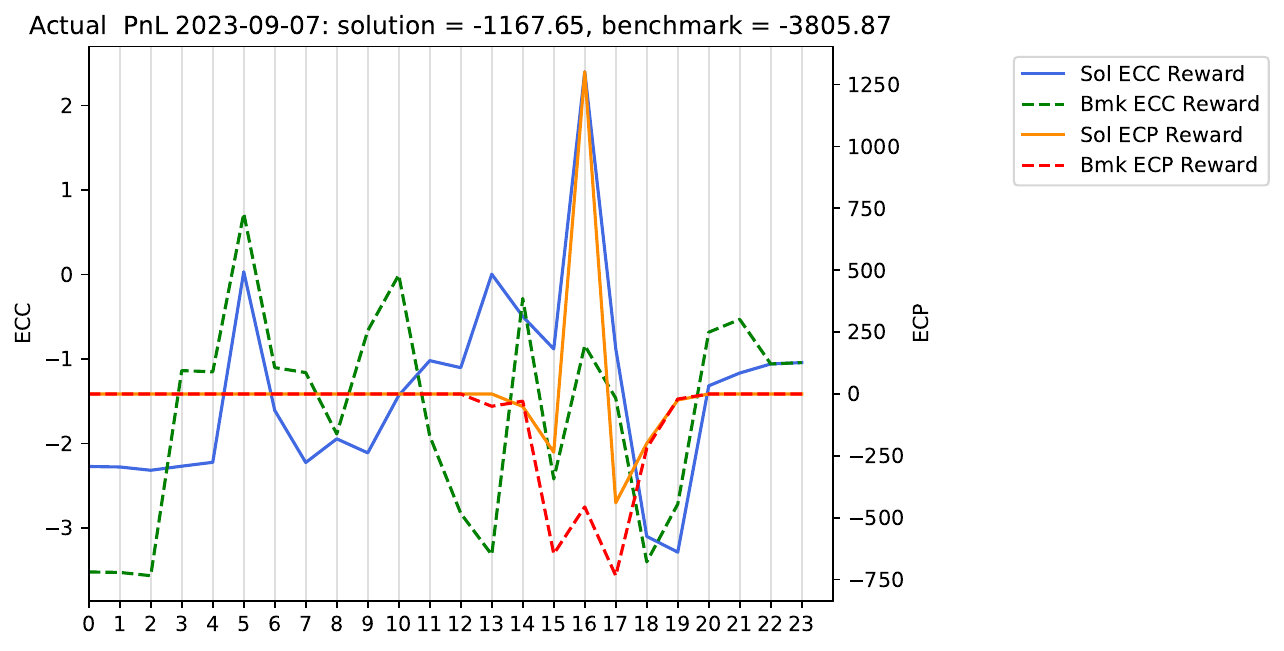}
    \caption{Performance on the final CP day, \texttt{2023-09-07}} 
    \label{sfig:high_cp_prob}
    \end{subfigure}
    \caption{Daily performance for a subset of days with the charging schedule of the BESS and the associated SOC (on the left), and the hourly rewards associated to the electricity consumption (ECC) and the coincident peak charge (ECP) (on the right). Note that the non-coincident peak charge is not represented but is included in the total actual P\&L for the solution and the benchmark.}
    \label{fig:full_cp_pnl_daily}
\end{figure}

For each of these days, Figure \ref{fig:full_cp_pnl_daily} represents the charging/discharging schedule of the battery and the associated State of Charge, as well as the hourly break down of the expected rewards from the electricity consumption and the coincident peak. This Profit and Loss function is computed at the end of the day $d$, using the knowledge of the microgrid's actual load $\hat{L}^{\texttt{MG}} = ( \hat{L}^{\texttt{MG}}_0, \dots, \hat{L}^{\texttt{MG}}_{23})$, the actual PV production $\widehat{PV} = (\widehat{PV}_0, \dots, \widehat{PV}_{23})$ and the NCP and CP probabilities from the method presented in \autoref{ssec:building_blocks}:

\begin{align}
    \operatorname{Actual \, P\&L}(d) = & \sum_{h=0}^{23} P_h \times \Big(\hat{L}^{\texttt{MG}}_{h} - \widehat{PV}_h - C (\pi^+_h - \pi^-_h) \Big) \nonumber  \\ \nonumber 
    & +  \lambda_{CP} \cdot p^{\texttt{PS}}_{d,nrm}\sum_{h=0}^{23} p^{\texttt{PS}}_h \Big( \hat{L}^{\texttt{MG}}_{h} - \widehat{PV}_h - C(\pi^+_h - \pi^-_h)\Big) \nonumber  \\ 
    &+ \lambda_{NCP} \cdot p^{\texttt{MG}}_{d,nrm} \Big( \hat{L}^{\texttt{MG}}_{h^*} -\widehat{PV}_{h^*} - C(\pi^+_{h^*} - \pi^-_{h^*}) \Big)
\end{align}

where $h^*$ is the hour of the maximum net load on the grid for that day, i.e.
$$
h^*=\operatorname{argmax}_{h=0,\cdots,23} \hat{L}^{\texttt{MG}}_{h} - \widehat{PV}_{h} - C(\pi^+_h - \pi^-_h).
$$

\autoref{sfig:low_prob} shows the behavior of both strategies for a day that has a zero coincident peak probability. Since there is no contribution from the CP-hour probabilities to the cost function, the whole time interval from 9AM to 5PM is used for charging after the first peak shave. As a result, both algorithms do not generate any coincident peak reward. We can however note that our solution provides a smoother electricity consumption, which is a desirable feature in terms of safety and robustness for the grid. \autoref{sfig:med_prob} shows a day with medium CP-day probability of $0.504$: the charging after the first peak shave is thus reduced to non-CP potential hours, i.e. hours that have a zero CP-hour probability. Since the potential CP reward is significant, the proposed algorithm will discharge the battery during the two hours of highest CP probability. The benchmark does not predict any CP-event and therefore does not dispatch the battery. 
\autoref{sfig:high_prob} shows a day where both methods predict a CP-event: the accounting of the CP rewards leads to a preference for the proposed solution. 
\autoref{sfig:high_cp_prob} shows the performance on the final CP-day, which is the true single CP of the whole period.  The benchmark does lead to a reduction of the CP cost, since it discharged the battery during the highest CP probability hours, but this still resulted in a cost to the microgrid. On the contrary, the solution leads to a significant reward due to a larger discharge during the hour of highest CP-probability. 

\begin{figure}[htb!]
    \centering
     \includegraphics[width=0.8\linewidth]{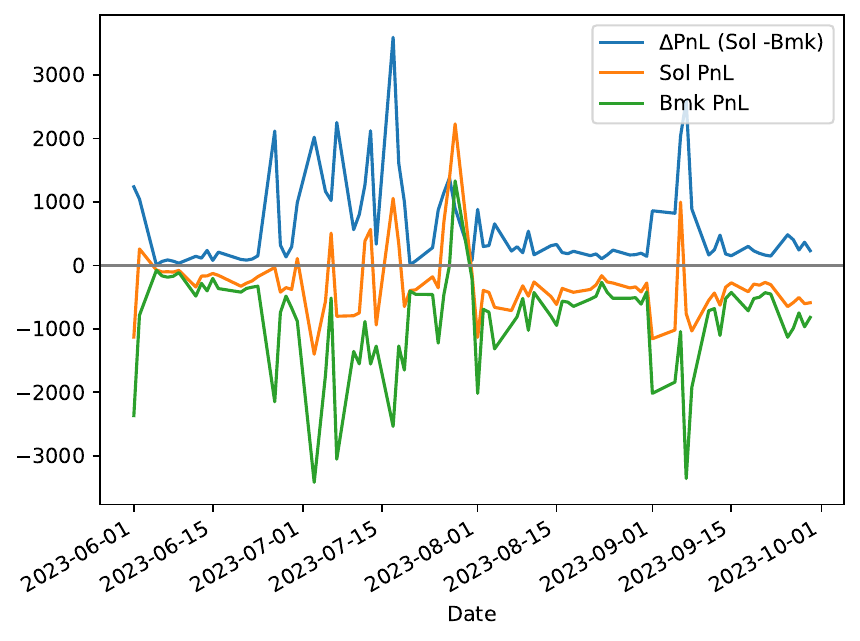}      
     \caption{Comparison of daily PnL without degradation}
     \label{fig:pnl_year_2023}
\end{figure}

\autoref{fig:pnl_year_2023} shows the timeseries of the total actual P\&L including Electricity Consumption, Coincident-Peak and Non-Coincident Peak charges, without the degradation loss. We notice that the proposed solution systematically does better than the benchmark, since the orange curve is always above the green one. The spikes in the blue curve, representing the difference of the solution and the benchamrk P\&L, are mostly associated with the CP reward collection.

\subsection{Final performance P\&L}

In this subsection, we reevaluate the daily P\&L after the summer period is over, i.e. once the single coincident peak of PSE\&G and the four non-coincident peaks can be determined. This function is computed sequentially using the exact knowledge of the running CP updates and the running NCP updates within each month. Hence, for a given day $d$, we consider: 

\begin{align}
    \operatorname{Final \, P\&L}(d) = & \sum_{h=0}^{23} P_h \times \Big(\hat{L}^{\texttt{MG}}_{h} - \widehat{PV}_h - C (\pi^+_h - \pi^-_h) \Big) \nonumber  \\ \nonumber 
    & +  \lambda_{CP} \mathbbm{1}_{\{d,CP\}}\cdot \sum_{h=0}^{23} \mathbbm{1}_{\{h,CP\}}\Big( \hat{L}^{\texttt{MG}}_{h} - \widehat{PV}_h - C(\pi^+_h - \pi^-_h)\Big) \nonumber  \\ 
    &+ \lambda_{NCP} \mathbbm{1}_{\{d,NCP\}} \cdot  \Big( \hat{L}^{\texttt{MG}}_{h^*} -\widehat{PV}_{h^*} - C(\pi^+_{h^*} - \pi^-_{h^*}) \Big)
\end{align}
where $\mathbbm{1}_{\{d,CP\}}$ (resp. $\mathbbm{1}_{\{d,NCP\}}$) is an indicator function taking the value $1$ if day $d$ sees a CP (resp. NCP) update and $0$ otherwise, $\mathbbm{1}_{\{h,CP\}}$ is similarly defined for the CP-hour. 

\autoref{tab:true_pnl_table} presents the break down of the P\&L according to the Electricity Consumption Cost, the demand charges and the transmission charge. We notice that the proposed algorithm leads to a significant reduction of each of the monthly demand charges, with an approximate change of 45\%. The most striking difference relates to the transmission charge, with the solution leading to more than \$4,000 in rewards against a cost of \$931 for the benchmark. 

\begin{table}[htb!]
    \centering
    \resizebox{\textwidth}{!}{
    \begin{tabular}{|c|c|c|c|} \hline 
     Algorithm & Benchmark & Solution & $\Delta$ (Solution - Benchmark) \\ \hline 
     Total Electricity Consumption Cost (USD) & 3,250.35 & 3,253.28 & 2.93 \\
     Demand charge for June  (USD)  & 3,477.81 & 1,875.90 & -1,601.91\\
     Demand charge for July  (USD) & 3,437.79 & 1,868.09 &  -1,569.70\\
     Demand charge for August   & 3,345.13 & 1,844.79 & -1,500.65\\
     Demand charge for September   & 3,509.21 & 1,928.70 & -1,580.51\\
     Final Coincident Peak Charge  & 931.76 & - 4,116.17 & -5,029.93 \\ \hline
    Total Final Cost & 17,952.05 & 6,654.59 &  -11,294.49\\ \hline 
    \end{tabular}
    }
    \caption{Final Profit and Loss (P\&L) in USD computed at the end of the summer period of 2023 knowing the true single coincident peak of PSE\&G and the true monthly non-coincident peaks of the microgrid. }
    \label{tab:true_pnl_table}
\end{table}

\begin{figure}[htb!]
    \centering 
     \begin{subfigure}[t]{0.45\textwidth}
        \includegraphics[width=\linewidth,height=5cm]{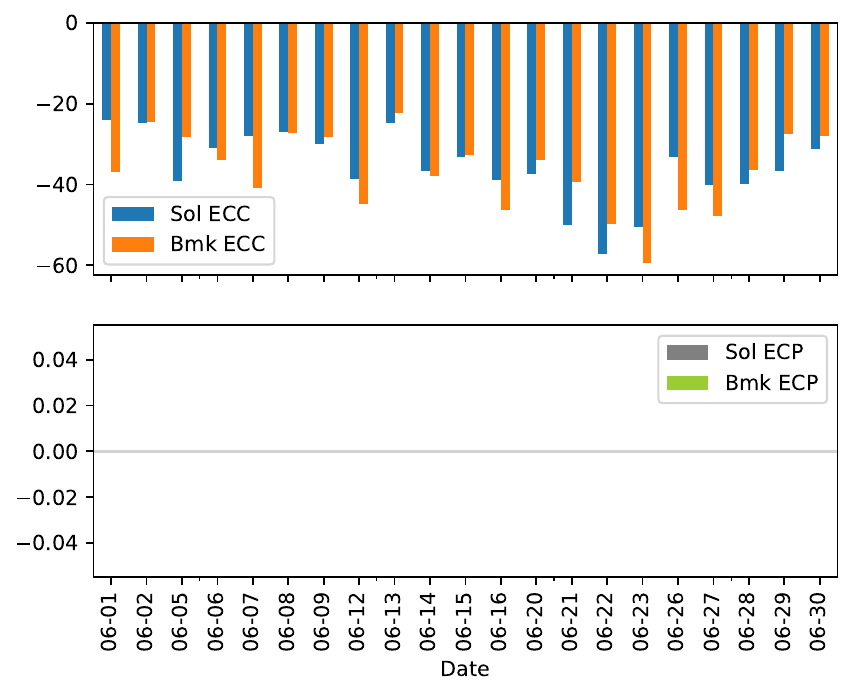}
    \end{subfigure}
     \begin{subfigure}[t]{0.45\textwidth}
     \includegraphics[width=\linewidth,height=5cm]{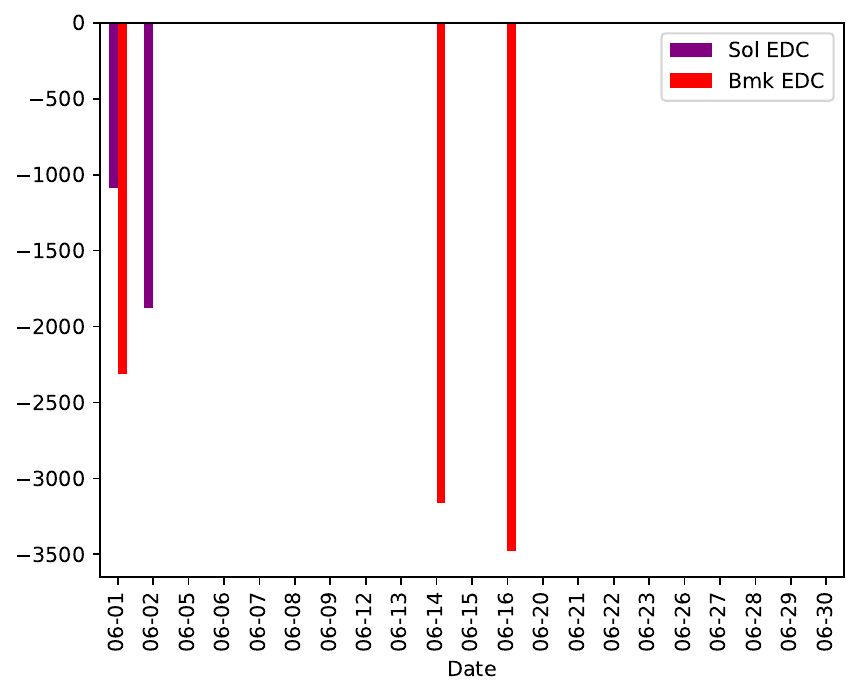}
     \end{subfigure}
    
    \begin{subfigure}[t]{0.45\textwidth}
        \includegraphics[width=\linewidth,height=5cm]{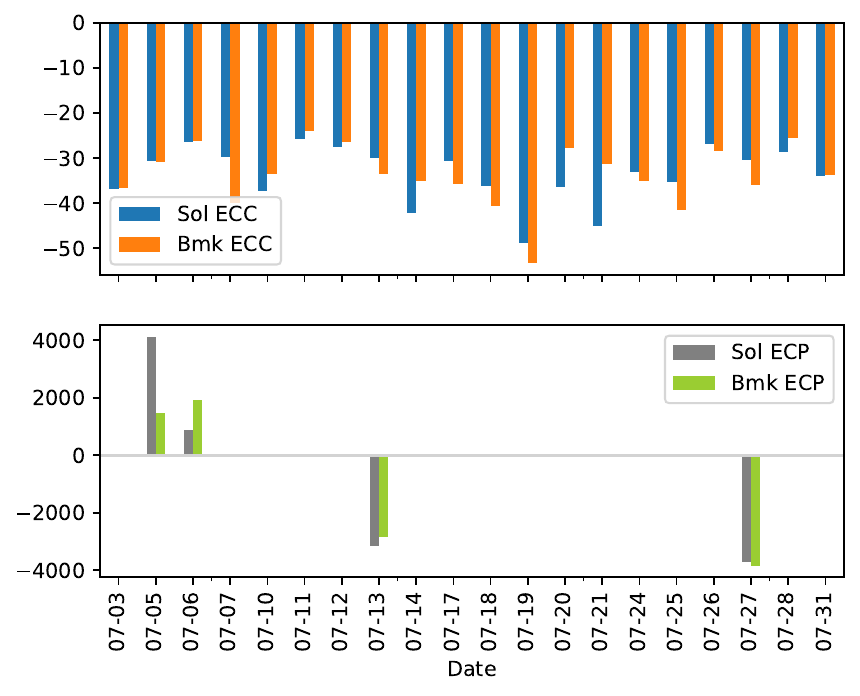}
    \end{subfigure}
     \begin{subfigure}[t]{0.45\textwidth}
     \includegraphics[width=\linewidth,height=5cm]{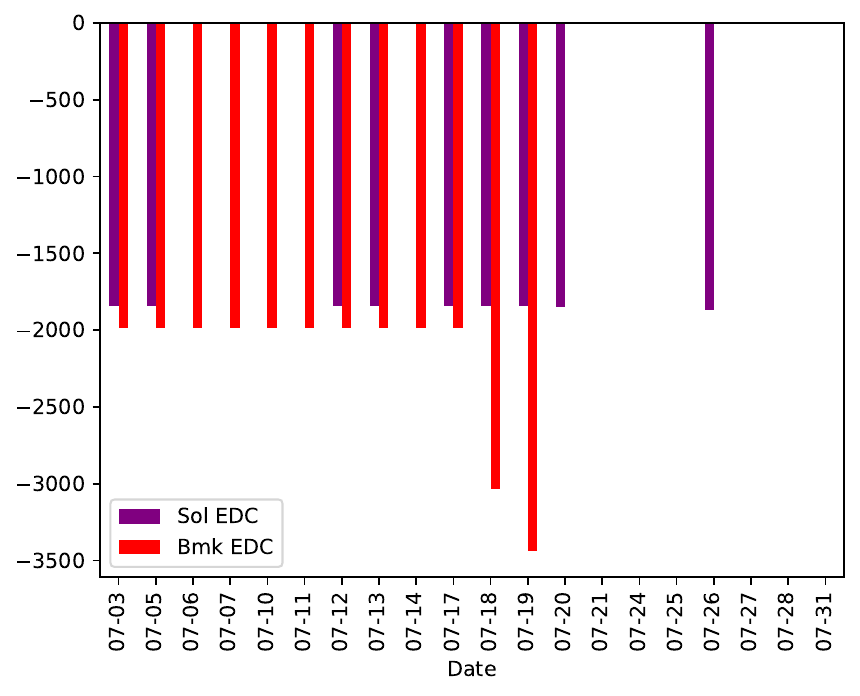}
     \end{subfigure}

          \begin{subfigure}[t]{0.45\textwidth}
        \includegraphics[width=\linewidth,height=5cm]{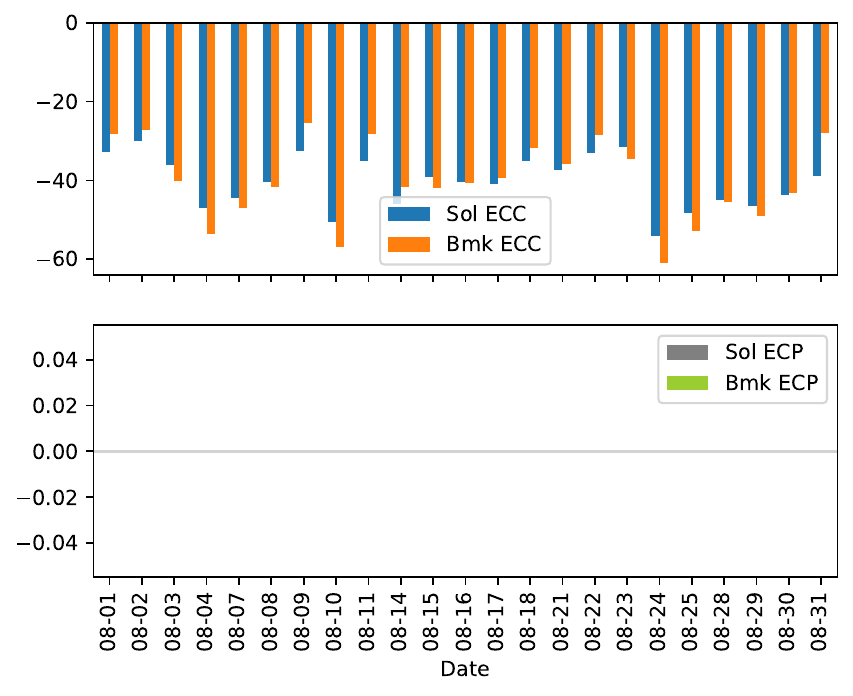}
    \end{subfigure}
     \begin{subfigure}[t]{0.45\textwidth}
     \includegraphics[width=\linewidth,height=5cm]{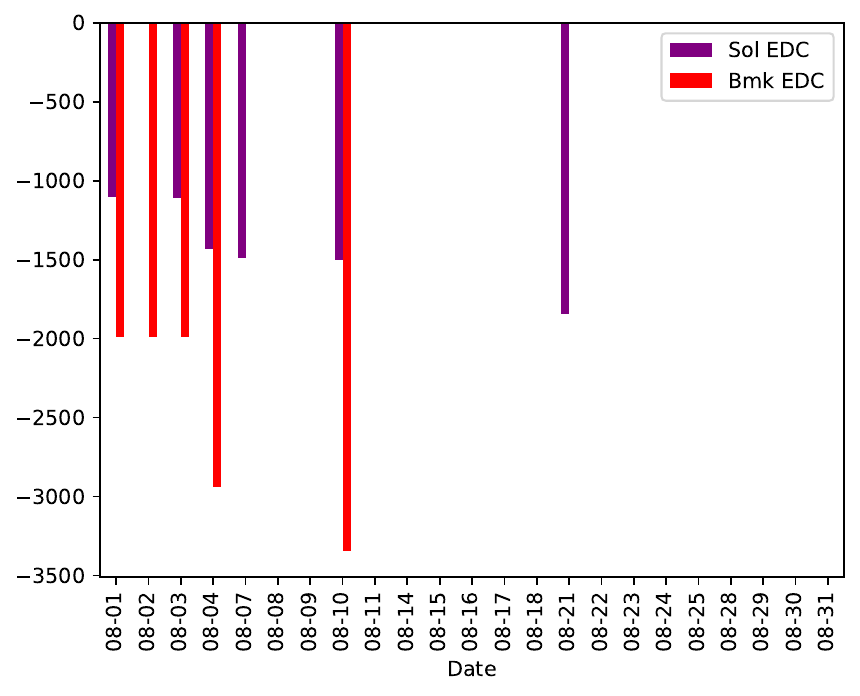}
     \end{subfigure}
          \begin{subfigure}[t]{0.45\textwidth}
        \includegraphics[width=\linewidth,height=5cm]{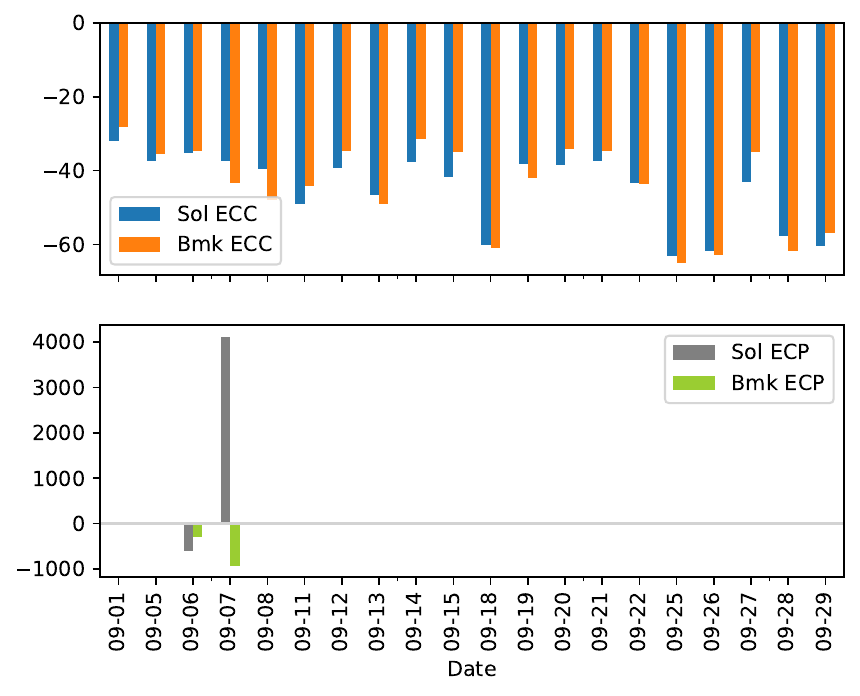}
    \end{subfigure}
     \begin{subfigure}[t]{0.45\textwidth}
     \includegraphics[width=\linewidth,height=5cm]{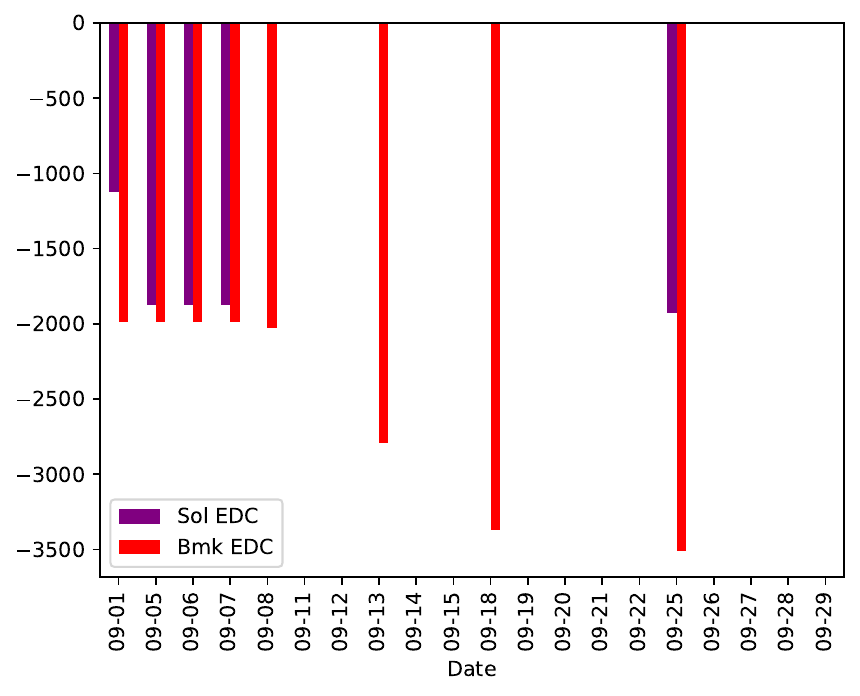}
     \end{subfigure}
     \caption{True PnL broken down by month for 2023: the difference of the electricity consumption cost (ECC) and that of the expected coincident peak charge (ECP) between the solution and the benchmark  are illustrated day by day on the left, and the demand charges (EDP) for each Non-Coincident Peak update are illustrated on the right. }
     \label{fig:true_pnl_diff}
\end{figure}

\autoref{fig:true_pnl_diff} shows the daily variations of these bill items. The histograms with orange and blue bars show the Electricity Consumption Reward: this consumption is driven by similar trends, with differences bounded by 15 in absolute value. The histograms with gray and green bars show the difference of CP-based reward for the real CP-updates, namely whenever a new running maximum is seen. The final CP charge corresponds to the last bars, which illustrates the CP reward earned through the proposed solution against the resulting cost for the benchmark. Finally, the bar plots on the right show the updates of the demand rewards each time a new running maximum for the load (including the additional consumption for the battery) is seen for the current month. The most notable difference is that the proposed solution leads to an effective reduction of non-coincident peaks with a stabilization of their value, compared to the increasing costs of the benchmark. This shows that the proposed methodology effectively reduces the demand charges. 

\section{Conclusion} 
In this work, we focused on a battery management optimization problem to increase the reward of a CP program and simultaneously, minimize the NCP demand charges. This is the more important that both can represent a significant portion of electricity bills for commercial and industrial customers. Our approach requires a simulation engine capable of generating Monte-Carlo scenarios to predict coincident peaks based on the ISO load, and non-coincident peaks, based on the local load. Based on these scenarios, we were able to formulate a mixed integer non-linear program for the battery dispatch schedule, which incorporated the probability of occurrence of coincident peak and non coincident peak events. The superior performance of our method was demonstrated against a sophisticated benchmark relying on historical load profile for a microgrid located in the jurisdiction of PSE\&G and PJM; its superiority to shave peaks and reduce transmission charges was demonstrated. Provided that data are available, our method can be easily adapted to all utility networks of PJM, and other power grids. 

\section*{Acknowledgments}
Authors R.C., X.Y. and C.Z. were partially supported by ARPA-
E grant DE-AR0001289. They would like to thank Xiaofan Wu from Siemens Technology for insightful discussions about coincident peak programs. 

\printbibliography

\end{document}